\catcode`\@=11

\font\tenmsa=msam10 \font\sevenmsa=msam7 \font\fivemsa=msam5
\font\tenmsb=msbm10
\font\sevenmsb=msbm7 \font\fivemsb=msbm5 \newfam\msafam \newfam\msbfam
\textfont\msafam=\tenmsa \scriptfont\msafam=\sevenmsa
\scriptscriptfont\msafam=\fivemsa \textfont\msbfam=\tenmsb
\scriptfont\msbfam=\sevenmsb \scriptscriptfont\msbfam=\fivemsb

\def\hexnumber@#1{\ifnum#1<10 \number#1\else \ifnum#1=10 A\else\ifnum#1=11
 B\else\ifnum#1=12 C\else \ifnum#1=13 D\else\ifnum#1=14 E\else\ifnum#1=15
 F\fi\fi\fi\fi\fi\fi\fi}

\def\msa@{\hexnumber@\msafam} \def\msb@{\hexnumber@\msbfam}
\mathchardef\boxdot="2\msa@00 \mathchardef\boxplus="2\msa@01
\mathchardef\boxtimes="2\msa@02 \mathchardef\square="0\msa@03
\mathchardef\blacksquare="0\msa@04 \mathchardef\centerdot="2\msa@05
\mathchardef\lozenge="0\msa@06 \mathchardef\blacklozenge="0\msa@07
\mathchardef\circlearrowright="3\msa@08 \mathchardef\circlearrowleft="3\msa@09
\mathchardef\rightleftharpoons="3\msa@0A
\mathchardef\leftrightharpoons="3\msa@0B \mathchardef\boxminus="2\msa@0C
\mathchardef\Vdash="3\msa@0D \mathchardef\Vvdash="3\msa@0E
\mathchardef\vDash="3\msa@0F \mathchardef\twoheadrightarrow="3\msa@10
\mathchardef\twoheadleftarrow="3\msa@11 \mathchardef\leftleftarrows="3\msa@12
\mathchardef\rightrightarrows="3\msa@13 \mathchardef\upuparrows="3\msa@14
\mathchardef\downdownarrows="3\msa@15 \mathchardef\upharpoonright="3\msa@16
 \mathchardef\downharpoonright="3\msa@17
\mathchardef\upharpoonleft="3\msa@18 \mathchardef\downharpoonleft="3\msa@19
\mathchardef\rightarrowtail="3\msa@1A \mathchardef\leftarrowtail="3\msa@1B
\mathchardef\leftrightarrows="3\msa@1C \mathchardef\rightleftarrows="3\msa@1D
\mathchardef\Lsh="3\msa@1E \mathchardef\Rsh="3\msa@1F
\mathchardef\rightsquigarrow="3\msa@20
\mathchardef\leftrightsquigarrow="3\msa@21 \mathchardef\looparrowleft="3\msa@22
\mathchardef\looparrowright="3\msa@23 \mathchardef\circeq="3\msa@24
\mathchardef\succsim="3\msa@25 \mathchardef\gtrsim="3\msa@26
\mathchardef\gtrapprox="3\msa@27 \mathchardef\multimap="3\msa@28
\mathchardef\therefore="3\msa@29 \mathchardef\because="3\msa@2A
\mathchardef\doteqdot="3\msa@2B 
\mathchardef\traceiangleq="3\msa@2C \mathchardef\precsim="3\msa@2D
\mathchardef\lesssim="3\msa@2E \mathchardef\lessapprox="3\msa@2F
\mathchardef\eqslantless="3\msa@30 \mathchardef\eqslantgtr="3\msa@31
\mathchardef\curlyeqprec="3\msa@32 \mathchardef\curlyeqsucc="3\msa@33
\mathchardef\preccurlyeq="3\msa@34 \mathchardef\leqq="3\msa@35
\mathchardef\leqslant="3\msa@36 \mathchardef\lessgtr="3\msa@37
\mathchardef\backprime="0\msa@38 \mathchardef\risingdotseq="3\msa@3A
\mathchardef\fallingdotseq="3\msa@3B \mathchardef\succcurlyeq="3\msa@3C
\mathchardef\geqq="3\msa@3D \mathchardef\geqslant="3\msa@3E
\mathchardef\gtrless="3\msa@3F \mathchardef\sqsubset="3\msa@40
\mathchardef\sqsupset="3\msa@41
\mathchardef\trianglelefteq="3\msa@45 \mathchardef\bigstar="0\msa@46
\mathchardef\between="3\msa@47 \mathchardef\blacktriangledown="0\msa@48
\mathchardef\blacktriangleright="3\msa@49
\mathchardef\blacktriangleleft="3\msa@4A
\mathchardef\blacktriangle="0\msa@4E \mathchardef\triangledown="0\msa@4F
\mathchardef\eqcirc="3\msa@50 \mathchardef\lesseqgtr="3\msa@51
\mathchardef\gtreqless="3\msa@52 \mathchardef\lesseqqgtr="3\msa@53
\mathchardef\gtreqqless="3\msa@54 \mathchardef\Rrightarrow="3\msa@56
\mathchardef\Lleftarrow="3\msa@57 \mathchardef\veebar="2\msa@59
\mathchardef\barwedge="2\msa@5A \mathchardef\doublebarwedge="2\msa@5B
\mathchardef\angle="0\msa@5C \mathchardef\measuredangle="0\msa@5D
\mathchardef\sphericalangle="0\msa@5E \mathchardef\varpropto="3\msa@5F
\mathchardef\smallsmile="3\msa@60 \mathchardef\smallfrown="3\msa@61
\mathchardef\Subset="3\msa@62 \mathchardef\Supset="3\msa@63
\mathchardef\Cup="2\msa@64  \mathchardef\Cap="2\msa@65
 \mathchardef\curlywedge="2\msa@66
\mathchardef\curlyvee="2\msa@67 \mathchardef\leftthreetimes="2\msa@68
\mathchardef\rightthreetimes="2\msa@69 \mathchardef\subseteqq="3\msa@6A
\mathchardef\supseteqq="3\msa@6B \mathchardef\bumpeq="3\msa@6C
\mathchardef\Bumpeq="3\msa@6D \mathchardef\lll="3\msa@6E 
\mathchardef\ggg="3\msa@6F  \mathchardef\circledS="0\msa@73
\mathchardef\pitchfork="3\msa@74 \mathchardef\dotplus="2\msa@75
\mathchardef\backsim="3\msa@76 \mathchardef\backsimeq="3\msa@77
\mathchardef\complement="0\msa@7B \mathchardef\intercal="2\msa@7C
\mathchardef\circledcirc="2\msa@7D \mathchardef\circledast="2\msa@7E
\mathchardef\circleddash="2\msa@7F \def\ulcorner{\delimiter"4\msa@70\msa@70 }
\def\urcorner{\delimiter"5\msa@71\msa@71 }
\def\llcorner{\delimiter"4\msa@78\msa@78 }
\def\lrcorner{\delimiter"5\msa@79\msa@79 } \def\yen{\mathhexbox\msa@55 }
\def\checkmark{\mathhexbox\msa@58 } \def\circledR{\mathhexbox\msa@72 }
\def\maltese{\mathhexbox\msa@7A } \mathchardef\lvertneqq="3\msb@00
\mathchardef\gvertneqq="3\msb@01 \mathchardef\nleq="3\msb@02
\mathchardef\ngeq="3\msb@03 \mathchardef\nless="3\msb@04
\mathchardef\ngtr="3\msb@05 \mathchardef\nprec="3\msb@06
\mathchardef\nsucc="3\msb@07 \mathchardef\lneqq="3\msb@08
\mathchardef\gneqq="3\msb@09 \mathchardef\nleqslant="3\msb@0A
\mathchardef\ngeqslant="3\msb@0B \mathchardef\lneq="3\msb@0C
\mathchardef\gneq="3\msb@0D \mathchardef\npreceq="3\msb@0E
\mathchardef\nsucceq="3\msb@0F \mathchardef\precnsim="3\msb@10
\mathchardef\succnsim="3\msb@11 \mathchardef\lnsim="3\msb@12
\mathchardef\gnsim="3\msb@13 \mathchardef\nleqq="3\msb@14
\mathchardef\ngeqq="3\msb@15 \mathchardef\precneqq="3\msb@16
\mathchardef\succneqq="3\msb@17 \mathchardef\precnapprox="3\msb@18
\mathchardef\succnapprox="3\msb@19 \mathchardef\lnapprox="3\msb@1A
\mathchardef\gnapprox="3\msb@1B \mathchardef\nsim="3\msb@1C
\mathchardef\napprox="3\msb@1D
\mathchardef\nsupseteqq="3\msb@23 \mathchardef\subsetneqq="3\msb@24
\mathchardef\supsetneqq="3\msb@25
\mathchardef\supsetneq="3\msb@29 \mathchardef\nsubseteq="3\msb@2A
\mathchardef\nsupseteq="3\msb@2B \mathchardef\nparallel="3\msb@2C
\mathchardef\nmid="3\msb@2D \mathchardef\nshortmid="3\msb@2E
\mathchardef\nshortparallel="3\msb@2F \mathchardef\nvdash="3\msb@30
\mathchardef\nVdash="3\msb@31 \mathchardef\nvDash="3\msb@32
\mathchardef\nVDash="3\msb@33 \mathchardef\ntrianglerighteq="3\msb@34
\mathchardef\ntrianglelefteq="3\msb@35 \mathchardef\ntriangleleft="3\msb@36
\mathchardef\ntriangleright="3\msb@37 \mathchardef\nleftarrow="3\msb@38
\mathchardef\nrightarrow="3\msb@39 \mathchardef\nLeftarrow="3\msb@3A
\mathchardef\nRightarrow="3\msb@3B \mathchardef\nLeftrightarrow="3\msb@3C
\mathchardef\nleftrightarrow="3\msb@3D \mathchardef\divideontimes="2\msb@3E
\mathchardef\varnothing="0\msb@3F \mathchardef\nexists="0\msb@40
\mathchardef\mho="0\msb@66 \mathchardef\thorn="0\msb@67
\mathchardef\beth="0\msb@69 \mathchardef\gimel="0\msb@6A
\mathchardef\daleth="0\msb@6B \mathchardef\lessdot="3\msb@6C
\mathchardef\gtrdot="3\msb@6D \mathchardef\ltimes="2\msb@6E
\mathchardef\rtimes="2\msb@6F \mathchardef\shortmid="3\msb@70
\mathchardef\shortparallel="3\msb@71 \mathchardef\smallsetminus="2\msb@72
\mathchardef\thicksim="3\msb@73 \mathchardef\thickapprox="3\msb@74
\mathchardef\approxeq="3\msb@75 \mathchardef\succapprox="3\msb@76
\mathchardef\precapprox="3\msb@77 \mathchardef\curvearrowleft="3\msb@78
\mathchardef\curvearrowright="3\msb@79 \mathchardef\digamma="0\msb@7A
\mathchardef\varkappa="0\msb@7B \mathchardef\hslash="0\msb@7D
\mathchardef\hbar="0\msb@7E \mathchardef\backepsilon="3\msb@7F
\def\Bbb{\ifmmode\let\next\Bbb@\else
\def\next{\errmessage{Use \string\Bbb\space only in math mode}}\fi\next}
\def\Bbb@#1{{\Bbb@@{#1}}} \def\Bbb@@#1{\fam\msbfam#1}

\catcode`\@=\active

 \def\CC{\hbox{{$\cal C$}}}

\def\CR{\hbox{{$\cal R$}}} 
 \def\CV{\hbox{{$\cal V$}}}
\def\CM{\hbox{{$\cal M$}}}

\def\cu{\hbox{{\sl u}}} 
\def\cv{\hbox{{\sl v}}}

\def\Vec{{\rm Vec}}

\def\lform{\hbox{$\sqcup$}\llap{\hbox{$\sqcap$}}}

\def\h{{{1\over2}}}

\def\A{{\Bbb A}}

\def\C{{\Bbb C}}
\def\Z{{\Bbb Z}}

\def\eps{{\epsilon}}

\def\aut{{\rm Aut\, }}

\def\rcross{{\triangleright\!\!\!<}}
\def\lcross{{>\!\!\!\triangleleft}}
\def\rcocross{{\blacktriangleright\!\!<}}
\def\lcocross{{>\!\!\blacktriangleleft}}

\def\rbiprod{{\cdot\kern-.33em\triangleright\!\!\!<}}
\def\lbiprod{{>\!\!\!\triangleleft\kern-.33em\cdot\, }}

\def\tens{\mathop{\otimes}}

\def\la{{\triangleright}}\def\ra{{\triangleleft}}

\def\isom{{\cong}}

\def\Aut{{\rm Aut}}

\def\Nat{{\rm Nat}}
\def\Ad{{\rm Ad}}

\def\id{{\rm id}}

\def\<{\langle}
\def\>{\rangle}

\def\equad{\kern -1.7em}

\def\eqn#1#2{\begin{equation}#2\label{#1}\end{equation}}

\def\o{{}_{\scriptscriptstyle(1)}}
\def\t{{}_{\scriptscriptstyle(2)}}
\def\th{{}_{\scriptscriptstyle(3)}}
\def\fo{{}_{\scriptscriptstyle(4)}}

\def\bo{{}^{\bar{\scriptscriptstyle(1)}}}
\def\bt{{}^{\bar{\scriptscriptstyle(2)}}}
\def\Ro{{\CR^{\scriptscriptstyle(1)}}}
\def\Rt{{\CR^{\scriptscriptstyle(2)}}}

\def\und#1{{\underline {#1}}}

\def\uo{{{}^{\scriptscriptstyle(1)}}}
\def\ut{{{}^{\scriptscriptstyle(2)}}}

\def\Bo{{{}_{\und{\scriptscriptstyle(1)}}}}
\def\Bt{{{}_{\und{\scriptscriptstyle(2)}}}}

\def\text#1{\mbox{\rm #1}}
\def\note#1{}

\def\blacksquare{{\lform}}
\def\frac#1#2{{{#1\over#2}}}

\def\proof{\goodbreak\noindent{\bf Proof\quad}}

\def\endproof{{\ $\lform$}\bigskip }

\def\alignn#1#2{\begin{eqnarray}\label{#1}#2
\end{eqnarray}}
\def\cmath#1{\[\begin{array}{c} #1 \end{array}\]}
\def\ceqn#1#2{\begin{equation}\label{#1}
\begin{array}{c}#2\end{array}\end{equation}}

\def\vecu{{\bf u}}\def\vecx{{\bf x}}

\documentstyle[11pt, epsf]{article}
\newtheorem{lemma}{Lemma}[section] \newtheorem{propos}[lemma]{Proposition}
\newtheorem{example}[lemma]{Example} \newtheorem{theorem}[lemma]{Theorem}

\begin{document}\baselineskip 20pt

{\ }\qquad\qquad \hskip 4.3in 
\vspace{.2in}

\begin{center} {\LARGE  SOME COMMENTS ON BOSONISATION AND BIPRODUCTS}
\\ \baselineskip 13pt{\ }
{\ }\\
S. Majid\footnote{Royal Society University Research Fellow and Fellow of
Pembroke College, Cambridge}\\
{\ }\\
Department of Mathematics, Harvard University\\
Science Center, Cambridge MA 02138, USA\footnote{During the calendar years 1995
+ 1996}\\
+\\
Department of Applied Mathematics \& Theoretical Physics\\
University of Cambridge, Cambridge CB3 9EW\\
\end{center}

\vspace{20pt}
\begin{quote}\baselineskip 13pt
\noindent{\bf Abstract}
We collect here some less well-known results and formulae about the
bosonisation construction which turns braided groups into quantum groups. We
clarify the relation with biproduct Hopf algebras (the constructions are not
the same), the response to twisting of braided groups and the abstract
characterisation  via automorphisms of the forgetful functor for the category
of (co)modules of a braided group.

\bigskip
\noindent Keywords:  quantum group -- braided group -- twisting -- bosonisation
-- biproduct -- crossed module -- braided reconstruction -- colour Lie algebra
-- supersymmetry

\end{quote}
\baselineskip 19.5pt

\section{Introduction}

Recently there has been some interest in the theory of braided groups or Hopf
algebras in braided categories\cite{Ma:bra}\cite{Ma:bg} and the bosonisation
functor which relates them to quantum groups\cite{Ma:bos}. Applications in
physics include the spectrum generating quantum groups\cite{MacMa:spe}, the
construction of inhomogeneous quantum groups\cite{Ma:poi} and the cross product
structure of the quantum double\cite{Ma:dou}\cite{Ma:skl}\cite{Ma:mec}.
Applications in pure mathematics include~\cite{CFW:sch} and~\cite{FisMon:sch}.
Here we collect some modest results about the construction and address some
frequently asked questions. As a novel feature, we give all four possible
versions (left modules, right modules, left comodules, right comodules) of the
various formulae. This should make the paper quite useful as a reference. We
also correct a mathematical confusion in the recent J. Algebra
paper~\cite{FisMon:sch} where it was assumed incorrectly that bosonisation and
the related theory of biproducts\cite{Rad:str} have the same input data.  We
provide a natural counterexample to this assertion. The main new result of the
paper is a detailed calculation of the automorphism braided group
$BGL_q(2)\rcocross \A_q^2$  of the forgetful functor from the category of
braided $\A_q^2$, which demonstrates explicitly a step from the abstract
construction of bosonisaton~\cite{Ma:bos} (in comodule form). This construction
fully solves a question recently posed by B.~Pareigis\cite{Par:rec} about
`hidden symmetries'.

We work over a ground field $k$ and use the usual notations $S,\eps$ for the
antipode and counit, and  $\Delta h= h\o\tens h\t$ for the coproduct applied to
an element $h\in H$\cite{Swe:hop} (summations understood). We use the symbols
$\lcross$, etc., for cross (or smash) products, $\lcocross$ for cross (or
smash) coproducts and $\lbiprod$ when both are made simultaneously.

This is the final version of a preprint with similar title and the same
mathematical content, circulated in April 1995. Only some material about
twisting in Section~3 and assossiativity (\ref{crossassoc}) has been added.

\section{Module and Comodule Formulae}

One of the joys of Hopf algebras is that for every theorem of a certain general
type one gets three theorems free. We will develop this as a formal result in
Section~5, but for the moment we merely demonstrate the principle at work for
braided groups. Thus, in the original work on braided groups\cite{Ma:bra} we
worked with braided groups $B$ living in the category of left-modules of a
quasitriangular Hopf algebra $H,\CR$. We used this version of the theory
because it is more familiar for physicists.  Less well-known perhaps is  the
dual version which was also introduced by the author\cite{Ma:bg}, in which we
work with $B$ in the category of right-comodules of a dual-quasitriangular Hopf
algebra $H$. Even less well-known is a theory for right modules or left
comodules.

By quasitriangular bialgebra or Hopf algebra we mean a Hopf algebra $H$
equipped with invertible $\CR\in H\tens H$ obeying the axioms of
Drinfeld\cite{Dri}
\ceqn{qua}{(\Delta\tens\id)\CR=\CR_{13}\CR_{23},\quad
(\id\tens\Delta)\CR=\CR_{13}\CR_{12}\\ \Delta^{\rm op}=\CR(\Delta\ )\CR^{-1}.}
If we consider $\CR$ as a map $k\to H\tens H$ and write Drinfeld's axioms as
commuting diagrams, and then reverse all arrows, we have the dual concept of a
dual-quasitriangular (or coquasitriangular) Hopf
algebra\cite{Ma:rep}\cite[Thm.~4.1]{Ma:pro}. One can then write those
axioms out explicitly
as a skew bialgebra bicharacter $\CR:H\tens H\to k$ with respect to which $H$
is quasi-commutative (the dual of Drinfeld's axioms). Explicitly:
\ceqn{dqua}{\CR(hg\tens f)=\CR(h\tens f\o)\CR(g\tens f\t),\quad \CR(h\tens
gf)=\CR(h\o\tens f)\CR(h\t\tens g)\\
g\o h\o\CR(h\t\tens g\t)=\CR(h\o\tens g\o)h\t g\t.}
Drinfeld requires $\CR$ invertible, so we also have to dualise this concept as
$\CR$ invertible in the convolution algebra of maps $H\tens H\to k$.

The rest of Drinfeld's theory can also be dualised. If some results are
routinely formulated with diagrams (such as cross product constructions) this
is just a matter of reversing arrows. Or if they involve complicated algebra it
 may be easier to prove the dual version directly. For example, among less
well-known results about dual-quasitriangular Hopf algebras one finds in an
appendix to~\cite{Ma:bg}:

\begin{propos}\cite[Prop~A.5]{Ma:bg} Let $H$ be a dual quasitriangular Hopf
algebra. Then the square of the antipode is inner in the convolution algebra
$H\to H$ and hence the antipode is bijective.
\end{propos}

It means that the assumption that the antipode of a dual quasitriangular Hopf
algebra $H$ is bijective, which is usually assumed by pure mathematicians
e.g.~\cite[Thm~2.15]{FisMon:sch}, should be deleted as superfluous.

Likewise, if the category of modules of a quasitiangular Hopf algebra is
braided, then so is the category of comodules of a dual-quasitriangular one.
Similarly, this time reflecting all diagrams about a vertical axis, we have
both left and right versions of the theory. We recall that a braided category
means a braiding $\Psi$ between any two objects. We denote left actions by
$\la$ and right actions by $\ra$. We denote coactions by $v\mapsto v\bo\tens
v\bt$, which lives in $H\tens V$ for a left comodule $V$ and $V\tens H$ for a
right comodule. So the braiding in the four cases is:
\ceqn{psi}{ \Psi_L(v\tens w)=\Rt{}\la w\tens \Ro{}\la v,\quad \Psi_R(v\tens
w)=w\ra\Ro\tens v\ra\Rt\\
\Psi^L(v\tens w)=\CR(w\bo\tens v\bo)w\bt\tens v\bo,\quad \Psi^R(v\tens
w)=w\bo\tens v\bo\CR(v\bt\tens w\bt).}
The four categories are denoted ${}_H\CM$, $\CM_H$, ${}^H\CM$ and $\CM^H$, for
left, right modules and left, right comodules respectively. This is basically
in Drinfeld\cite{Dri} and more explicitly in \cite[Sec. 7]{Ma:qua} in the
module version, among other works from about this time.

Braided groups make sense as algebraic structures within any braided category,
but these four categories are the most important. Indeed, super and colour Lie
algebras and Hopf algebras have all been studied in isolation for many years.
One of the main ideas of the theory of braided groups is (as well as to
generalise them to the braided case) to cast these concepts  as constructions
in one of these four categories. Thus, in~\cite{Ma:tra}\cite{Ma:exa} we
introduced the quantum group $\Z_2'$ which generates the category of
supervector spaces as ${\rm SuperVec}={}_{\Z_2'}\CM$. This was generalised to
the category of anyonic (or $\Z_n$-graded) vector spaces in
1991\cite{Ma:any} and further on to categories generated by Abelian groups
equipped with bicharacters\cite{Ma:csta}, which is the setting for colour Lie
algebras, except that we do not assume that the bicharacter is skew (in the
skew case the category is symmetric rather than braided and the theory is more
straightforward). While one can work with these categories directly,
appreciating that they are generated by a quantum group allows one to apply
some general Hopf algebra machinery, such as the bosonisation
theorem\cite{Ma:bos}. A common misconception is that this use of
quantum groups as generating categories in which we can do such algebra was
developed previously in the symmetric, e.g. super case, before the advent of
braided groups. As far as I know, it is due to the author and useful,
e.g.\cite{MacMa:spe}, even in the super case.

Next we turn to the concept of braided group itself, as well as other algebraic
constructions in braided categories. A braided group $B$ is like a Hopf algebra
except that the coproduct $\und\Delta:B\to B\und\tens B$ is a homomorphism to
the {\em braided tensor product} algebra. The product here involves the
braiding $\Psi$. In concrete cases it is
\eqn{btens}{ (b\tens c)(a\tens d)=b\Psi(c\tens a)d.}
In our four preferred categories, the homomorphism property becomes:
\ceqn{copmult}{\und\Delta(bc)=b\Bo (\CR\ut\la c\Bo)\tens (\CR\uo\la
b\Bt)c\Bt,\quad \und\Delta(bc)=b\Bo ( c\Bo\ra\CR\uo)\tens (
b\Bt\ra\CR\ut)c\Bt\\
\und\Delta(bc)= \CR(c\Bo\bo\tens b\Bt\bo)b\Bo c\Bo\bt\tens
b\Bt\bt c\Bt,\quad \und\Delta(bc)= b\Bo c\Bo\bo\tens b\Bt\bo c\Bt
\CR(b\Bt\bt\tens c\Bo\bt),}
where $\und\Delta b=b\Bo\tens b\Bt$ denotes the braided coproduct. Similarly,
the antipode for braided groups is a braided-antihomomorphism:

\begin{propos}\cite[Fig.~2]{Ma:tra} The antipode $S$ of a braided group $B$
is a braided-antihomorphism in the sense
\[ \und S\circ \cdot=\cdot\circ\Psi\circ(\und S\tens \und S),\quad
\und\Delta\circ
\und S=(\und S\tens \und S)\circ\Psi\circ\und\Delta.\]
\end{propos}

Less well-known, however, is that these properties are not the axioms of a
braided group but rather they require proof (one of the first non-trivial
lemmas in braided group theory).  In fact, the antipode axioms are  $(\und
Sb\Bo)b\Bt=\und\eps(b)=b\Bo \und S b\Bt$ as usual (said diagrammatically in a
general category). In our four preferred categories the braided
antimultiplicativity property in Proposition~2.2 becomes:
\ceqn{Sant}{ \und S(bc)=(\CR\ut\la \und Sc)(\CR\uo\la\und Sb),\quad \und
S(bc)=( \und Sc\ra\CR\uo)( \und Sb\ra\CR\ut)\\
\und S(bc)= \CR(c\bo\tens b\bo)(\und S c\bt)(\und S b\bt),\quad \und S(bc)=
(\und S c\bo)(\und S b\bo)\CR(b\bt\tens c\bt).}
For example, it is perhaps not clear in~\cite[eqn.~(1.9)]{FisMon:sch} that the
homomorphism property for $\und\Delta$ is part of the definition, while the
property for $\und S$ is not part of the definition but follows as explained
above.

The proof of Proposition~2.2 in~\cite{Ma:tra} is by diagrammatic means. As far
as I know, it is the only known proof, direct algebraic proofs being
impractical.
In this diagrammatic method one writes products as
$\cdot=\epsfbox{prodfrag.eps}$, coproducts and
$\und\Delta=\epsfbox{deltafrag.eps}$ and the braiding as
$\Psi=\epsfbox{braid.eps}$ with inverse $\Psi^{-1}=\epsfbox{braidinv.eps}$. In
this way algebraic information `flows' along braid and knot
diagrams\cite{Ma:bra}\cite{Ma:tra}\cite{Ma:bos}, not unlike the manner in which
information flows along the wires in a computer. Such wiring diagrams are a
standard feature in mathematics and engineering, and have even been used for
ordinary Hopf algebras\cite{Yet:rep}. In all these previous contexts there is
no non-trivial operator attatched to the crossing of wires. One just wires
outputs into inputs without caring about whether one passes over or under
another. The novel feature of braided groups is that now, for the first time,
crossings  represent non-trivial operators $\Psi$ or $\Psi^{-1}$. The notation
makes sense by combining standard ideas about wiring diagrams with the
coherence theorem for braided categories\cite{JoyStr:bra}. These $\Psi$
correspond to the quasitriangular structure (\ref{psi}), which is precisely the
key complication in braided group formulae such as (\ref{copmult}) and
(\ref{Sant}). This is why it was indispensable in
\cite{Ma:bra}\cite{Ma:tra}\cite{Ma:bos}.

In these terms we would like to address a further misconception that  the
braided group theory  follows  automatically as a generalisation of the   older
theory of Hopf algebras in symmetric categories. In fact, the symmetric theory
follows in a canonical way from the theory of ordinary Hopf algebras since one
merely inserts a `symmetry' $\Psi$ in place of transposition in every usual
algebraic construction, e.g.~\cite{Gur:yan}\cite{Par:non}. The braided case is
much more problemmatic because not only must one choose from $\Psi,\Psi^{-1}$
(they do not coincide in the braided case), but there {\em may be no consistent
choice at all}, i.e. a standard construction for quantum groups may simply get
`tangled up' in the braided setting. Even such basic things as the tensor
product of braided groups (within the catgeory) and the Jacobi identity in its
usual form become tangled up in the braided case; they fail and a
new theory is needed\cite{Ma:lie}. The diagrammatic notation is one of the main
tools introduced  in the braided theory to help control this problem. There is,
however, no automatic way to go from usual results about quantum groups or Hopf
algebras in symmetric categories to strictly braided ones. For this reason we
really should distinguish carefully between algebraic constructions in
symmetric categories, for which there is a canonical procedure to extend
general categorical constructions to this case, and the braided case which
requires genuinely new work.

Finally, we would like to clarify some notational confusion for which the
author is certainly to blame.  In the first `transmutation' construction
$B(\ ,\ )$ for braided groups we emphasised the diagonal case $B(H,H)$ which
are braided-cocommutative\cite{Ma:bra} or
braided-commutative\cite{Ma:exa}\cite{Ma:eul} in a certain sense, i.e. like
classical groups. These remain some of the most interesting  for conformal
field theory~\cite{Ma:rec}\cite{Lyu:mod}\cite{LyuMa:bra}, but they are only one
  case; we have subsequently used the term `braided groups' for any Hopf
algebra in a braided category (not only the strict usage as braided
(co)commutative ones). For example, the general $B(\ ,\ )$   transmuation yield
quantum-braided groups  with braided (dual)quasitriangular
structure\cite{Ma:bg}\cite{Ma:tra}. Thus,~\cite{Ma:exa} emphasised for
physicists the braided matrices $B(R)=B(A(R),A(R))$, while more recently the
general cases $B(R,Z)=B(A(R),A(Z))$ have proven  interesting as
well\cite{Hla:bra}\cite{Lu:bra}\cite{MaPla:any}. Strictly speaking $B(R)$ and
$B(R,Z)$ are obtained by transmutation\cite{Ma:bg} only for suitably nice cases
where $A(Z)$ can be replaced by an actual Hopf algebra; but once the formulae
are obtained under this assumption, they can all be checked directly  along the
lines explained in~\cite{Ma:exa} in the diagonal case and
{}~\cite{Hla:bra}\cite{Lu:bra} in the general case. There are now also linear
braided groups\cite{Ma:poi} which are not of this transmutation type at all.

In summary, the braided group theory is  different in a fundamental way from
the theory of Hopf algebras in symmetric categories (where $\Psi=\Psi^{-1}$).
The first examples were introduced (by the author) in both module and comodule
form. The module version of the theory\cite{Ma:rec}\cite{Ma:bra} was presented
to the mathematical physics community in Kyoto in May 1990 and in St Petersburg
in September 1990. The comodule version\cite{Ma:eul}\cite{Ma:bg}\cite{Ma:tan}
was presented to the Hopf algebra community at the AMS meeting in San Francisco
in January, 1991. More general examples of braided groups (not necessarily
associated to quantum groups at all) are provided by a more general
automorphism construction due to the author\cite{Ma:rec}\cite{Ma:bg} and
independently in the diagonal case to Lyubashenko\cite{Lyu:mod}.

\section{Bosonisation and Twisting}

In this section we recall two important theorems about braided groups, namely
the bosonisation construction itself\cite[Thm~4.1]{Ma:bos}  and
twisting\cite{Ma:euc}\cite{Ma:qsta}, giving them now in all of our four
categories.

Bosonisation generalises the Jordan-Wigner bosonisation transform for
$\Z_2$-graded systems
in physics, and associates to every Hopf algebra $B$ in the braided category of
representations of $H$ an equivalent ordinary Hopf algebra $B\lbiprod H$ (left
handed cases) or $H\rbiprod B$ (right handed cases). If $B\in
{}_H\CM$  then $B\lbiprod H$ is defined as the cross product by
the canonical action $\la$ of $H$ (by which $B$ is an object) and a coproduct
built from the quasitriangular structure $\CR$ and this action. Dually, if
$B\in{}^H\CM$ then $B\lbiprod H$ is defined as the cross coproduct by the
canonical coaction of $H$ (by which $B$ is an object) and a product built from
the dual quasistriangular structure $\CR$ and the coaction. The explicit
formulae in our four categories are:
\ceqn{bos}{{}_H\CM:\quad hb=(h\o\la b)h\t,\quad \Delta b= b\Bo\CR\ut\tens
\CR\uo\la b\Bt,\quad Sb=(\cu\CR\uo\la \und Sb)S\CR\ut\\
\CM_H:\quad bh=h\o (b\ra h\t),\quad \Delta b= b\Bo\ra \CR\uo\tens \CR\ut
b\Bt,\quad Sb=(S\CR\ut)\und S b\ra\CR\uo\cv\\
{}^H\CM:\quad \Delta b=b\Bo b\Bt\bo\tens b\Bt\bt,\quad bh=\CR(h\o\tens
b\bo)b\bt h\t,\quad Sb=b\bo (\und S b\bt)\\
\CM^H:\quad\Delta b=b\Bo\bo\tens b\Bo\bt b\Bt,\quad bh=h\o b\bo\CR(b\bt\tens
h\t),\quad Sb=(\und S b\bo)b\bt}
In all cases, $H$ is a sub-Hopf algebra and $B$ a subalgebra. For the antipode,
$\cu=(S\CR\ut)\CR\uo$ and $\cv=\CR\uo(S\CR\ut)$.

There are several ways of thinking about this construction. The abstract
characterisation~\cite{Ma:bos} of the resulting Hopf algebras is that their
(co)modules are monoidally equivalent to the braided $B$-(co)modules in the
braided category. The module case is given explicitly in~\cite{Ma:bos}. The
dual theorem, which we will need in Section~6, is:

\begin{propos}cf\cite[Thm~4.2]{Ma:bos} Let $H$ be a dual quasitriangular Hopf
algebra and $B$ a Hopf algebra in $\CM^H$. The $B$-comodules in $\CM^H$ can be
identified canonically with $H\rbiprod B$-comodules as monoidal categories over
${}_k\CM$.
\end{propos}
\proof Cf\cite{Ma:bos} a $B$-comodule in the category means a vector space
which is an $H$-comodule and a $B$-comodule which intertwines the $H$-coaction.
The corresponding coaction of $H\rbiprod B$ consists of the $B$-coaction
followed by the $H$-coaction. Using the properties of a dual quasitriangular
structure one finds that this is an identification of monoidal categories (i.e.
that the tensor product of comodules is respected). \endproof

In more concrete terms, there are at least two concrete points of view on the
same formulae. The natural one in the context of (\ref{bos}) is that the
coproduct in the left module case is the braided tensor product coalgebra
$B\und\tens H_L$  where $H_L$ denotes the left regular representation.
Similarly, in the right comodule case the product is the braided tensor product
algebra $H^R\und\tens B$ as in (\ref{btens}), where $H^R$ denotes the right
regular corepresentation ($H$ as a right comodule by its coproduct). Similarly
for right modules and left comodules.

A second point of view on the same formulae is that the coproduct in the left
module case has a cross coproduct form by coaction $\beta(b)=\CR\ut\tens
\CR\uo\la b$, which is the {induced coaction} introduced earlier
in~\cite{Ma:dou} -- so the bosonisation is manifestly both a cross product and
cross coproduct\cite{Ma:skl}. Similarly, in the right comodule bosonisation
cases, the product is a cross product by the {induced action}  $b\ra h=
b\bo\CR(b\bt\tens h)$ via the dual-quasitriangular structure, etc. This point
of view connects with a more general biproduct construction (see Section~4). It
is {\em not} however, the point of view which captures the key properties of
bosonisatons. Nor is it the point of view by which the construction was first
introduced.

A second important theorem for braided groups is the theory of `gauge
equivalence' or twisting of braided groups. Such twisting for  Hopf algebras
was introduced in the work of Drinfeld\cite{Dri} and used in~\cite{GurMa:bra};
it was extended to braided groups in~\cite{Ma:euc}\cite{Ma:qsta}. Its
importance for physics is that many systems can appear algebraically different
but should really be equivalent in a physical sense. For example, there are
curently two `twistor' formulations of $q$-spacetime based on $2\times 2$
matrices, described by the algebras $R_{21}\vecx_1\vecx_2=\vecx_2\vecx_1R$ and
$R_{21}\vecu_1R\vecu_2=\vecu_2 R_{21}\vecu_1R$,
respectively~\cite{Ma:euc}\cite{Ma:exa}. The use of such a form as spacetime is
due to the author (and works for a general R-matrix), but the choice of $su_2$
R-matrix recovers previous algebras proposed for `Euclidean'  and   `Minkowski'
$q$-coordinates (and now as braided groups). The use of this R-matrix form for
Minkowski space explicitly occured in 1992 in~\cite{Ma:mec}, with its braided
coaddition structure (also in R-matrix form) appearing in 1993 in Meyer's
paper~\cite{Mey:new}. The point of our discussion is that the twisting theory
of braided groups exactly relates the two systems; they are gauge equivalent at
the algebraic level, differing up to equivalence only in their choice of
$*$-structure. This makes possible the concept of `quantum Wick rotation'
between the two systems\cite{Ma:euc}. There are  many other applications of
twisting besides this one.

In the module setting, the data for twisting is $\chi\in H\tens H$ which is a
2-cocycle for $H$ in the sense
$\chi_{23}(\id\tens\Delta)\chi=\chi_{12}(\Delta\tens\id)\chi$ and
$(\eps\tens\id)\chi=1$ and ensures that $H_\chi$ with coproduct
$\Delta_\chi=\chi(\Delta\ )\chi^{-1}$ and quasitriangular structure
$\CR_\chi=\chi_{21}\CR\chi^{-1}$ is also a quantum group, the twisting of $H$.
In the comodule setting the dual data is $\chi:H\tens H\to k$ obeying
$\chi(h\o\tens f\o)\chi(g\tens h\t f\t)=\chi(g\o\tens h\o)\chi(g\t h\t\tens f)$
and $\chi(1\tens h)=\eps(h)$ and ensures that $H_\chi$ with product
$h\cdot_\chi g=\chi(h\o\tens g\o)h\t g\t \chi^{-1}(h\th\tens g\th)$ and dual
quasitriangular structure ${\CR}_\chi(h\tens g)=\chi(g\o\tens h\o)\CR(h\t\tens
g\t) \chi^{-1}(h\th\tens g\th)$ is also a quantum group, the dual-twisting of
$H$. More details are in~\cite{Ma:book}. Now, if $B$ is a braided group in one
of our preferred quantum-group generated braided categories then its twisting
$B_\chi$ lives in the category generated by the twisted quantum group. In the
twisting of braided groups, both the product and coproduct are modified. The
formulae in the four cases are
\ceqn{Btwist}{{}_{H_\chi}\CM:\quad  b\cdot_\chi c=\cdot\left( \chi^{-1}\la
(b\tens c)\right),\quad \und\Delta_\chi=\chi\la\und\Delta\\
 \CM_{H_\chi}:\quad b\cdot_\chi c=\cdot\left( (b\tens
c)\ra\chi^{-1}\right),\quad \und\Delta_\chi=(\und\Delta\ )\ra\chi\\
{}^{H_\chi}\CM:\quad b\cdot_\chi c= \chi^{-1}(b\bo\tens c\bo) b\bt c\bt , \quad
\und\Delta_\chi b=\chi(b\Bo\bo\tens b\Bt\bo) b\Bo\bt\tens b\Bt\bt\\
\CM^{H_\chi}:\quad b\cdot_\chi c=  b\bo c\bo \chi^{-1}(b\bt\tens c\bt), \quad
\und\Delta_\chi b=b\Bo\bo\tens b\Bt\bo\chi(b\Bo\bt\tens b\Bt\bt).}
The braided antipode, unit and counit do not change. The twisting formulae for
braided groups  appeared in~\cite{Ma:euc}, with the detailed proof that the
result is again a braided group appearing in~\cite{Ma:qsta}.

The twisting of braided groups commutes with their bosonisation. Thus
$B_\chi\lbiprod H_\chi\isom (B\lbiprod H)_\chi$ and $H_\chi\rbiprod B_\chi\isom
(H\rbiprod B)_\chi$ are also (long) theorems proven in~\cite{Ma:qsta}. Here
$\chi$ is viewed in the bosonised algebra in the trivial way in order to make
the twisting after bosonisation.

Here we want to mention a possible application of these ideas to colour Hopf
algebras and Lie algebras\cite{Sch:gen}. These can be understood as algebraic
structures in the comodule category generated by $kG,\beta$ where
$\beta:G\times G\to k$ is a bicharacter extended as a dual-quasitriangular
structure. Here $G$ is Abelian, and the case usually studied is when
$\beta(s,t)=\beta(t,s)^{-1}$ (the skew case). In this case the resulting
category (of $G$ graded spaces with transposition defined by $\beta$) is
symmetric and life is much easier. The point is that some colour Lie algebras,
while appearing genuinely different, may be twisting equivalent to usual ones.
If so then their algebraic properties will tend to be equivalent as well; one
can prove results about them by twisting to the `gauge' where they become
ordinary Lie algberas,  using a theorem about them there, and twisting back to
obtain the corresponding theorem for the original colour Lie algebra.  We have:

\begin{propos} If a skew bicharacter $\beta=\chi^2$ for some other skew
bicharacter $\chi$, then a colour Lie algebra in the category generated by
$kG,\beta$ is twisting equivalent to a usual Lie algebra.
\end{propos}
\proof We twist $kG,\beta$ by $\chi$. Then
$\beta_\chi(s,t)=\chi(t,s)\beta(s,t)\chi^{-1}(s,t)=\chi(t,s)\chi(s,t)=1$. The
category generated by $kG,\beta$ therefore twists to the category of $G$-graded
vector spaces with its usual trivial transposition. All consructions in the
category twist. We consider the enveloping colour Hopf algebra\cite{Gur:yan}.
It has coproduct $\und\Delta\xi=\xi\tens1+1\tens\xi$ and relations
$\xi\eta-\beta(|\xi|,|\eta|)\eta\xi=[\xi,\eta]$ for homogeneous elements of
degree $|\ |$ of the colour Lie algebra (working here with right $kG$
comodules, say). Twisting gives us the relations
$\xi\cdot \eta -\eta\cdot\xi=\chi^{-1}(|\xi|,|\eta|)\xi\eta
-\chi(|\xi|,|\eta|)\eta\xi=\chi^{-1}(|\xi|,|\eta|)
(\xi\eta-\beta(|\xi|,|\eta|)\eta\xi)=\chi^{-1}(|\xi|,|\eta|)[\xi,\eta]\equiv
[\xi,\eta]_\chi$. It is easy to check that if $[\ ,\ ]$ obeys the colour Jacobi
identity etc. (defined in the obvious way with transposition $\beta$) then $[\
,\ ]_\chi$ obeys the usual one.  \endproof

More generally, if $\beta=\beta_0\chi^2$ for some other skew bicharacters
$\beta_0,\chi$ then the same calcultion shows that a $\beta$-colour Lie algebra
is twisting equiavalent to a $\beta_0$-colour Lie algebra by the same formulae.
The simplest among skew bicharacters are those which have values $\pm1$, which
makes them both skew and symmetric. They are super-like in the sense that the
transposition is generalised by $\pm1$ factors.  Many   colour Lie algebras
can be `reduced' up to twisting to ones of this super-like type.

\begin{example} Let $G=(\Z/m\Z)^n$ and $\beta(s,t)=q^{(s,t)}$ for some
antisymmetric $\Z/m\Z$-valued bilinear form on $G$ and $q^m=1$, a primitive
$m$-th root of 1. A $\beta$-colour Lie algebra is twisting equivalent to an
ordinary Lie algebra if $m\equiv 1,3$ mod 4. It is twisting equivalent to a
super-like Lie algebra if $m\equiv 2$ mod 4.
\end{example}
\proof If $m$ is odd then $2$ is invertible in $\Z/m\Z$. Hence we can write
$\chi(s,t)=q^{\h(s,t)}$ and have an example of the preceding propositon. If
$m\equiv 2$ mod 4, we can write every element of $\Z/m\Z$ uniquely in the form
$i\equiv 2j$ or $i\equiv 2j+\h m$, where $j\in\{0,1,\cdots, \h m-1\}$. We do
this for the values of the bilinear form on the standard basis elements in the
upper-triangular range (i.e. we write our bilinear form as an antisymmetric
matrix and consider its upper-triangular entries). This defines two
antisymmetric matrices  and hence two antisymmetric bilinear forms $(\ ,\ )_1$,
$(\ ,\ )_0$ such that $(s,t)=2(s,t)_1+(s,t)_0$, where $(s,t)_0$ has values in
$\{0,\h m\}$. Then $\beta=\chi^2\beta_0$ where $\beta_0=q^{(\ ,\ )_0}$ and
$\chi=q^{(\ ,\ )_1}$. Here $\beta_0$ has values in $\pm1$. \endproof

It seems likely that some ideas of Scheunert\cite{Sch:gen} about reducing
certain commutation factors to super ones could be formulated as a twisting
equivalence along similar lines. This remains, however, for further work.

The extension of colour Hopf algebras to the braided (non-skew case) occured
in~\cite{Ma:any}\cite{Ma:csta}. In the latter we studied the bosonisation of
such  braided groups as a novel approach to physical quantisation.
For example, the braided line bosonises to the quantum plane\cite{Ma:csta}. At
the Lie algebra level, the braided case is much more complicated but can be
developed in the framework of~\cite{Ma:lie}.

\section{Bosonisation and Biproducts are Not the Same}

Some years ago, Radford characterised Hopf algebras which are both a cross
product and cross coproduct by an action of a Hopf algebra $H$ (what he called
`biproducts') as Hopf algebras equipped with a split projection to
$H$\cite{Rad:str}. It was shown by the author in~\cite{Ma:dou}\cite{Ma:skl}
that the acted-upon object $B$ as in fact a braided group in the braided
category ${}_H^H\CM$ of crossed modules (when $H$ has bijective antipode). It
was also explained that when $H$ is finite dimensional this category is just
the braided category of modules ${}_{D(H)}\CM$, which was already known by
then. Here $D(H)$ is Drinfeld's quantum double quasitriangular Hopf algebra.
Finally, it was shown that when $H$ is quasitriangular then the bosonisation
construction can be viewed as an example of this more general construction by
means of a certain functor ${}_H\CM\to {}_H^H\CM$. These results are all due to
the author\cite{Ma:dou}\cite{Ma:skl}. In this section we show that this functor
is not, however, an isomorphism (so the constructions are not the same),
contrary to recent assumptions in the literature \cite{FisMon:sch}. The
category ${}^H_H\CM$ itself was introduced in a different context in
\cite{Yet:rep}.  It  is also an example of a more general construction of the
`Pontryagin dual'\cite{Ma:rep} or `double'  of any monoidal category.

Consider for the moment $H$ finite-dimensional. Drinfeld's quantum
double\cite{Dri} $D(H)$ is a Hopf algebra containing (in some conventions)
$H,H^{*\rm op}$. So ${}_{D(H)}\CM$ just means left $H$-modules and left
$H^{*\rm op}$-modules which are compatible in that they respect the cross
relations of the double. But a left $H^{*\rm op}$-module is trivially the same
thing (by evaluation) as a left $H$-comodule. This is the category ${}^H_H\CM$.
Because $D(H)$ is quasitriangular, we know that this category is braided when
$H$ has bijective antipode, as appreciated independently in \cite{Yet:rep}.
Similarly, the category $\CM_{D(H)}$ can be formulated as $\CM^H_H$, consisting
of compatible right $H$ action and coaction.  Explicitly, the left and right
compatibility conditions and braidings are:
\ceqn{crossmod}{{}^H_H\CM:\quad h\o v\bo\tens h\t\la v\bt= (h\o\la v)\bo
h\t\tens (h\o\la v)\bt,\quad \Psi_L(v\tens w)=v\bo\ra w\tens v\bt\\
\CM^H_H:\quad v\bo\ra h\o\tens v\bt h\t= (v\ra h\t)\bo\tens h\o(v\ra
h\t )\bt,\quad \Psi_R(v\tens w)=w\bo\tens v\ra w\bt.}
It should be  clear that since we have dispensed with $H^{*\rm op}$ itself  we
do not need to assume that $H$ is finite-dimensional. So, associated to any
Hopf algebra with bijective antipode one has these left and right crossed
module braided categories. We do not actually need the antipode of $H$ but only
the inverse or skew antipode, for $\Psi$ to be invertible. And we only require
the latter in order to have a standard braided-categorical setting (it is not
needed for the actual constructions below). Morphisms in all these categories
are maps which intertwine both the module and comodule structures.

Now consider $B$ a braided group in ${}^H_H\CM$. The conditions entailed in
this ensure that the cross product and cross coproduct $B\lbiprod H$
simultaneously by the action and coaction is an ordinary Hopf
algebra\cite{Ma:dou}\cite{Ma:skl}. Conversely, every Hopf algebra with split
projection to $H$ is isomorphic to one of the form $B\lbiprod H$ for $B$ a Hopf
algebra in ${}^H_H\CM$. This is the braided version of Radford's
theorem\cite{Ma:skl}. Similarly for $B\in \CM^H_H$ we have $H\rbiprod B$ and a
Hopf algebra with split projection to $H$ is also isomorphic to one of this
form.

Some authors have  wondered whether this braided version of Radford's theorem
adds to what was known in \cite{Rad:str}, apart from some terminology. Here we
would like to explain that the answer is affirmative. In fact, knowing that $B$
is a braided group carries much more information than the properties elucidated
in \cite{Rad:str}; it tells us that the operator controlling the exotic nature
of the algebra-coalgebra $B$ is a braiding $\Psi$, obeying the Yang-Baxter
relations. It tells us that the product and coproduct are well-behaved with
respect to this operator (functoriality of $\Psi$), and other key properies of
braided groups which are needed to prove such things as the
antimultiplicativity in Proposition~2.2 etc. None of these properties are
implicit or hinted at in \cite{Rad:str}. There are some examples of `exotic'
algebra-coalgebras $B$ in Radford's paper but without proving such properites
as braiding, they were not shown at that time to be braided groups. Finally, we
note that the braided version puts $B$ into a category with other objects in
it, allowing us to make categorical constructions involving $B$ and other
objects.

For example, if $B,C$ are algebras in ${}^H_H\CM$ then the braided tensor
product $B\und\tens C$ from (\ref{btens}) is again an algebra in ${}^H_H\CM$.
Explicitly, it has product   $(b\tens c)(a\tens d)=b(c\bo\la a)\tens c\bt d$,
as explained by the author in \cite[Prop. A.2]{Ma:skl}.
This is a generalisation of the usual concept of cross product because the  the
latter can be viewed as $B\lcross H=B\und\tens H_{\Ad}^L$, where   $H^L_{\Ad}$
is an algebra in ${}^H_H\CM$  by the adjoint action and left regular coaction.
This point of view has been used in the braided case (where $H$ itself is a
braided group) by Bespalov\cite{Bes:cro}. To see how such generalised cross
products easily arise, note that $B\lcross H\in {}^H_H\CM$ by the tensor
product action and coaction, because it is a braided tensor product in
${}^H_H\CM$. If we make a cross product again by this tensor product action of
$H$, associativity of braided tensor products tells us that
\eqn{crossassoc}{ (B\lcross H)\lcross H=(B\und\tens H)\und\tens H=B\und\tens
(H\und\tens H)=B\und\tens (H\lcross H)}
where  $H\lcross H$ is a cross product by the adjoint action. In other words,
usual cross products are not closed under associativity but they are when
viewed as more general braided tensor products. The biproduct  $B\lbiprod H$
itself is an algebra in ${}^H_H\CM$ according to this. It is also a coalgebra
in ${}^H_H\CM$  as the braided tensor coproduct $B\und\tens H^{\Ad}_L$.

We note that although \cite{Ma:dou}\cite{Ma:skl} emphasised working in the
category ${}_{D(H)}\CM$, assuming that $H$ is finite-dimensional, each of these
papers also explained at the relevant point how to proceed in the
infinite-dimensional case using ${}^H_H\CM$. Specifically, it was observed in
\cite{Ma:dou} (below Cor.~2.3)  that (\ref{crossmod}) was one of Radford's
principal conditions in \cite{Rad:str}. (The supplementary  conditions that $B$
is an $H$-module coalgebra and an
$H$-comodule algebra were ommitted in~\cite[Cor~2.3]{Ma:dou}  but should also
be understood.) And it was observed in \cite{Ma:skl} (in the proof of Prop~A.2)
that the second of Radford's principal conditions is the braided group
homomorphism property for the braiding in (\ref{crossmod}). Moreover, we would
like to say that it is not the case that \cite[Prop.~A.2]{Ma:skl} only asserts
the converse direction; the forward direction that $B\in {}^H_H\CM$ gives a
Hopf algebra $B\lbiprod H$ is an integral part of the  proposition and is
covered in the proof.

Another corollary of the braided version of Radford's construction was that it
exhibited more clearly the connection with bosonisation. This is provided by
functors of braided categories when $H$ is quasitriangular or dual
quasitriangular. For our four preferred categories the functors are
\ceqn{embed}{ {}_H\CM\hookrightarrow {}^H_H\CM,\quad (V,\la)\mapsto
(V,\la,\beta),\quad \beta(v)=\CR\ut\tens\CR\uo\la v\\
\CM_H\hookrightarrow \CM^H_H,\quad (V,\ra)\mapsto
(V,\ra,\beta),\quad \beta(v)=v\ra\CR\uo\tens\CR\ut\\
{}^H\CM\hookrightarrow {}^H_H\CM, \quad (V,\beta)\mapsto
(V,\beta,\la),\quad h\la v=\CR(v\bo\tens h)v\bt\\
\CM^H\hookrightarrow \CM^H_H, \quad (V,\beta)\mapsto
(V,\beta,\ra),\quad v\ra h=v\bo \CR(v\bt\tens h)}
The idea in each case is to start with an action or coaction and {\em induce}
from it a compatible coaction or action. As far as I know, the first functor
is due to the author in \cite[Prop.~3.1]{Ma:dou} (including the
infinite-dimensional case), with the others as right-module or comodule
versions of the same result.
Because these are functors of braided categories, a braided group $B\in
{}_H\CM$, say, can be viewed   in ${}^H_H\CM$.  It is clear that the
corresonding bosonisation can be viewed as an example of a biproduct from this
second point of view.

Recently, this biproduct point of view on bosonisation was  emphasised
in~\cite{FisMon:sch}, although attributing the bosonisation construction
entirely to Radford\cite{Rad:str}. Indeed, the authors assert throughout the
paper  \cite[p.594, eqn~(1.16), below Prop~1.15, Remark~1.16]{FisMon:sch} that
${}^H\CM={}^H_H\CM$ when $H$ is dual-quasitriangular, so that the constructions
are strictly identified (and due to Radford since his paper \cite{Rad:str}
occured some years earlier).  We refer to the introduction of \cite{FisMon:sch}
where the bosonisation  papers are not  mentioned at all.

We show now that this identification ${}^H\CM={}^H_H\CM$ can never hold unless
$H$ is trivial.

\begin{propos} Let $H$ be a quasitriangular Hopf algebra with $H\ne k$. Then
the functor ${}_H\CM\to {}^H_H\CM$ introduced in~\cite{Ma:dou} is never an
isomorphism. Likewise, let $H$ be a dual quasitriangular Hopf algebra with
$H\ne k$. Then ${}^H\CM\to {}^H_H\CM$ is never an isomorphism.
\end{propos}
\proof This is clear in the finite-dimensional case from the construction in
\cite{Ma:dou}, where this functor was introduced as pull back along a Hopf
algebra projection $D(H)\to H$. Since $D(H)$ as a vector space is $H^*\tens H$,
this can never be isomorphism. An isomorphism of categories would, by
Tannaka-Krein reconstruction, require such an isomorphism. This is the
conceptual reason. For a formal proof which includes the infinite-dimensional
case, consider $H\in {}^H_H\CM$ by the left regular coaction $\Delta$ and left
adjoint action. If in the image of the first functor (with $H$
quasitriangular), then $ h\o\tens h\t= \CR\ut \tens \CR\uo\o h
S\CR\uo\t$ for all $h$ in $H$. Applying $\eps$ to the second factor tells us
that $h=\eps(h)$ for all $h$, i.e. $H=k$. This object in ${}^H_H\CM$ {\em can}
be in the image of the second functor, but this is iff the dual quasitriangular
structure is trivial and $H$ commutative. On other hand, consider $H\in
{}^H_H\CM$ by the left regular action and left adjoint coaction. If in the
image of the second functor (with $H$ dual quasitriangular) then $hg= \CR(g\o
Sg\th\tens h)g\t$. Setting $g=1$ tells us that $h=\eps(h)$ again, hence $H=k$.
This object can be in the image of the first functor, but this is iff the
quasitriangular structure is trivial and $H$ cocommutative. \endproof

This means in turn that general `biproducts' associated to $B\in {}^H_H\CM$ are
much more general than the Hopf algebras obtained by bosonisation when $B\in
{}_H\CM$ for $H$ quasitriangular or $B\in {}^H\CM$ for $H$ dual
quasitriangular, the  bosonisatons having many key properties not holding for
general biproducts.

The author would like to note for the record that  \cite{Ma:skl} was circulated
in January~1992 and its original form is archived on ftp.kurims.kyoto-u.ac.jp
as kyoto-net/92-02-07-majid. The original version of \cite{FisMon:sch} appeared
somewhat later (date received June~1992) and was shown to the author in
July~1992 at the  L.M.S. Symposium on Non-Commutative Rings in Durham, England.
The simpler version  with ${}^H_H\CM$ appearing in
\cite[Prop.~1.15]{FisMon:sch} was explained by the author to Susan Montgomery
at this time as the correct formulation dual to the module bosonisation
\cite{Ma:bos}. We refer to \cite[Rem. 1.16]{FisMon:sch} where the original
June~1992 version is described. This is perhaps not clear from the published
\cite{FisMon:sch}.

\section{Dualisation as Convention}

The principle that certain types of contructions for Hopf algebras have dual
versions is clear but perhaps not as widely accepted as it should be. One often
finds dual versions of known theorems presented in the literature as new. In
this short section I would like to elevate this principle to  a mathematical
theorem. After this, it really should not be necessary to publish certain types
of theorems four times. I would like to argue in fact that just as it is
generally accepted that using a right-handed version of a left-handed result is
merely a matter of convention and does not entail a new theorem, so the
reversal of arrows in the dual formulation is likewise not more than a matter
of convention.

One of the common objections to this point of view from experts is that the
Hopf algebras of interest  may not be  finite-dimensional and so may not have
an appropriate dual Hopf algebra. This argument is based, however, on a
misconception: it is not any specific Hopf algebra which is being dualised but
the theorem itself; its axioms and proofs:

\begin{theorem} Let $T$ be a theorem whose premises, proofs and results are
expressed by commuting diagrams in the category $\Vec$ of vector spaces. Then
(i) $T^{\rm op}$ defined by reversing all arrows in $T$ is also a correct
theorem in the category $\Vec$. (ii) $\bar T$ defined by reflecting all
diagrams in a mirror is another correct theorem in $\Vec$.
\end{theorem}
\proof (i) Reverse all arrows. In categorical terms we make the construction in
the opposite category where arrows are reversed. If the theorem involves
assuming a certain element in $H$, consider it as a map $k\to H$. If a theorem
involves the transposition map, reverse it as the transposition map again. A
theorem involving an algebra becomes one involving a coalgebra. A theorem
involving an action becomes one involving a coaction, etc. The axioms of a Hopf
algebra are self-dual in this way. Thus a theorem involving a Hopf algebra and
an action becomes one involving a Hopf algebra and a coaction, etc. (ii) It is
assumed that tensor products in $\Vec$ are all written horizontally and the
reflection is in a vertical axis. In categorical terms we make the construction
in the category $\Vec$ equipped with the opposite tensor product. Thus, left
actions become right actions, etc. Again, the axioms of a Hopf algebra are
symmetric in this sense. \endproof

An obvious example is the theorem that a left $H$-module algebra $B$  leads to
an associative algebra $B\lcross H$, the cross product. The dual theorem is
that a left $H$-comodule coalgebra $C$ leads to an coassociative cross
coproduct $C\lcocross H$.

There is an obviouos generalisation  to theorems in other categories, e.g to
braided group constructions in braided categories. Since these are routinely
done in any case by a certain diagrammatic notation (see Section~2),
dualisation or left-right reflection is even more routine. In this diagrammatic
notation all morphisms are considered pointing generally downwards. The dual
and left-right reversed version of a braided group construction is given in
this case by simply turning the diagram proofs up-side-down.

The theorem does not mean that emphasising the dual version of a construction
is not useful for some application, but it is the application itself which
would be new. Note also that there can still be a problem if we want not the
dual theorem but the
actual dual algebra or coalgebra etc. to a given one, i.e. dual in the sense of
Hopf algebra duality (adjunction of maps in a rigid category). This can take
rather more work. For example, \cite{Ma:mec} contains the proof that if $A,H$
are dually paired quantum groups and $B,C$ dually paired braided groups (the
correct definition of the latter is not completely obvious, and is not
symmetric) then the comodule and module bosonisations $A\rbiprod B$ and
$C\lbiprod H$ are dually paired Hopf algebra. Another version with $A\rbiprod
C$ and $B\lbiprod H$ dually paired is in \cite{Ma:qsta}.

\section{Braided Reconstruction}

Apart from the four openning paragraphs, this section is the same as in the
version of the paper circulated in April 1995. We recall the more abstract
definition of the bosonisation theorem\cite{Ma:bos}, giving it explicitly in
the comodule form. We observe that this provides without any work a solution to
a problem recently posed in \cite{Par:rec}: we will see  that the `hidden
symmetry' coalgebra remarked in \cite[Cor~5.7]{Par:rec} is in fact a braided
group, with the structure of a certain braided group cross coproduct, and we
will compute a detailed example.

Let us recall that the first result in the theory of braided groups is
establishing their existence. While one can easily write down
axioms for them, the main problem, which was not solved until 1989, was
establishing that those axioms can be satisfied non-trivially.
In a symmetric category there is no problem since one can take the enveloping
algebra of a colour or other generalised Lie algebra, for
example\cite{Gur:yan}, but this was not at all possible in the braided case.
The problem was solved (by the author) by introducing  the {\em automorphism
braided group} $\Aut(\omega)$ of a monoidal functor
$\omega$\cite{Ma:rec}\cite{Ma:bg}. A particular categorical realisation
corresponding to the case $\omega=\id$ (more precisely, the identity inclusion
in a cocompletion) was also considered independently, in the following year, by
V.~Lyubashenko.  Using the automorphism braided group construction we were able
to not only prove existence but to actually compute concrete examples of
braided groups\cite{Ma:exa} by means of a construction which we called {\em
transmutation}. If $G\to H$ is a Hopf algebra map and $H$ is
dual-quasitriangular then $G$ transmutes to a braided group $B(G,H)\in \CM^H$,
obtained abstractly as the automorphism braided group of the push out functor
$\CM^G\to \CM^H$.

Note that we do not need here the most general version\cite{Ma:bg} of
$\aut(\omega)$ in which $\omega:\CC\to\CV$ is a monoidal functor between a
monoidal category $\CC$ and a braided one $\CV$. For existence of this braided
group one needs a representing object for the natural transformation
$\Nat(\omega,\omega\tens(\ ))$ which requires a degree of cocompleteness and
rigidity. In \cite{Ma:bg} we began for simplicity with the strongest assumption
that $\CV$ is rigid and cocomplete over $\CC$. Later on in the same paper
\cite[p.205]{Ma:bg} we dropped the rigidity assumption in favour of rigidity of
$\CC$. Another option is to require that the image of $\omega$ is rigid. Each
approach has some advantages. The construction\cite{Ma:bg} itself is
independent of these details and proceeeds as soon as certain {\em
representability conditions} are satisfied, in whatever way. It is this aspect
of the automorphism construction which is also the most useful in practice:
after obtaining the required formulae under the most convenient assumptions for
representability,  one can  verify directly the properties of $B(G,H)$ etc. by
usual algebraic methods\cite{Ma:tra}\cite[Appendix]{Ma:bg}.

Building on this work, \cite{Ma:bos} introduced the bosonisation procedure as a
kind of `adjoint' to transmutation. We have described it concretely in
Section~3, but its abstract characterisation is as follows. Proceeding in
comodule form: let $B$ be a braided group in $\CC=\CM^H$. It is known from
\cite{Ma:bg} that its category $\CC^B$ of braided $B$-comodules  in $\CC$ also
has a tensor product (is a monoidal category). By more usual Tannaka-Krein
arguments\cite{Sav:cat}\cite{Ulb:hop} one knows that this category is
equivalent to the the usual comodules over a usual Hopf algebra. This is the
characteristic property of bosonisation\cite[Thm.~4.2]{Ma:bos} as explained in
Proposition~3.1. In the four cases:
\eqn{bosrep}{ {}_B({}_H\CM)={}_{\textstyle B\lbiprod H}\CM,\quad
(\CM_H)_B=\CM_{\textstyle H\rbiprod B},\quad {}^B({}^H\CM)={}^{\textstyle
B\lbiprod H}\CM,
\quad (\CM^H)^B=\CM^{\textstyle H\rbiprod B}.}
The identifications are the obvious ones of the underlying vector spaces. In
order to obtain this result, we passed through a more conceptual argument using
cross products in the braided category and transmutation. Thus, both $B(H,H)$
and $B$ are braided groups and the former (say, in the comodule case) coacts on
$B$ by a coaction corresponding under transmutation to the coaction of $H$ on
$B$ as an object in $\CM^H$. Hence\cite{Ma:bos}  we can make a braided group
cross coproduct $B(H,H)\rcocross B$ with the braided tensor product algebra. We
then recoginise such braided group cross (co)products as the transmutation of
some ordinary Hopf algebras $H\rbiprod B$. In our four categories:
\ceqn{boscon}{ B(H,B\lbiprod H)=B\lcross B(H,H), \quad B(H,H\rbiprod
B)=B(H,H)\rcross B,\\  B(B\lbiprod H,H)=B\lcocross B(H,H),\quad
 B(H\rbiprod B,H)=B(H,H)\rcocross B,}
where the transmutation functor $B(\ ,\ )$ in each case is the one appropriate
to the category. This is the abstract construction\cite{Ma:bos} of bosonisation
and the reason that it has a cross product form combined with a braided tensor
coproduct in the module case, or cross coproduct form combined with a braided
tensor product in the comodule case, as we have seen in Section~3.

Now we want to observe that this construction solves  automatically the
question posed in \cite{Par:rec}, namely what braided group does one
reconstruct as $\Aut(\omega)$ when we are given $B\in\CM^H$ and the forgetful
functor $\omega:(\CM^H)^B\to \CM^H$?

\begin{propos} Let $H$ be dual quasitriangular and $B$ a Hopf algebra in
$\CM^H$. Then the forgetful functor $\omega$ from $B$-comodules in $\CM^H$ to
$\CM^H$ has as automorphisms the braided group $B(H,H)\rcocross B$
in $\CM^H$. It has the cross coproduct coalgebra and
braided tensor product algebra, and is a transmutation of the bosonisation
$H\rbiprod B$ of $B$.
\end{propos}
\proof  Under the equivalence (\ref{bosrep}), the forgetful functor $\omega$
becomes the functor induced by push-out along the canonical Hopf algebra map
$H\rbiprod B\to H$ defined by the counit of $B$. But the automorphism braided
group of a functor induced by push out is exactly the definition of the
transmutation construction $B(\ ,H)$. So the answer is exactly the
transmutation $B(H\rbiprod B,H)$. But the abstract definition of bosonisation
 in (\ref{boscon}) means that this is just
$B(H,H)\rcocross B$. Indeed, these are exactly the conceptual steps (in
comodule form)
which led to the bosonisation theory~\cite{Ma:bos} in the first place.
\endproof

This demonstrates how one may use bosonisation theory: we
convert our problem for the braided group $B$ to one for its equivalent
ordinary Hopf algebra $H\rbiprod B$. Explicitly, the braided group $B(H,H)$
associated to $H$ is obtained as the automorphism braided group of the
identity functor from $\CM^H$ to itself\cite{Ma:eul}\cite{Ma:bg} and
corresponds to $B=k$. Its structure is $H$ as a coalgebra, with the right
adjoint coaction and  modified product\cite{Ma:bg}
\eqn{BHH}{  h\bo\tens h\bt= h\t \tens (Sh\o)h\th,\quad h\cdot g=
h\t g\t \CR((Sh\o)h\th\tens Sg\o)}
in terms of the structure of $H$. We consider $B(H,H)$ coacting on any $B$ by
the same
map $\beta$ by which $H$ coacts on $B$ as an object (the {\em tautological
coaction}). The fact that one can then make a (braided) cross coproduct by this
 and still obtain a Hopf algebra in the braided category with the braided
tensor product algebra structure reflects the fact that $B(H,H)$ is
braided-commutative with respect to $B$ in a certain (unobvious) sense
introduced in\cite{Ma:bg}. This was the key idea behind the construction in
\cite[Sec.~2]{Ma:bos}. This $B(H,H)\rcocross B$ has product and coproduct
defined diagrammatically cf.\cite[Sec.~2]{Ma:bos}
\eqn{rcross}{\epsfbox{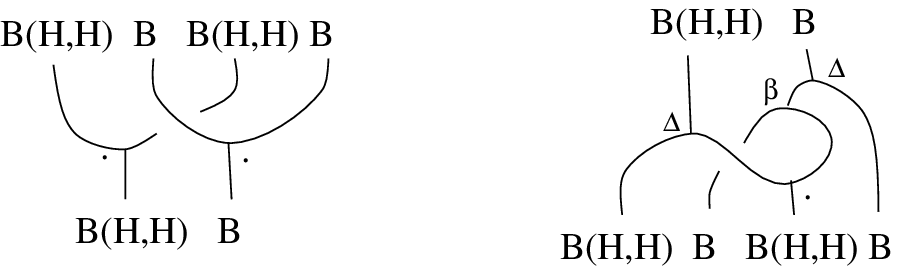}}
where $\Psi=\epsfbox{braid.eps}$ is the braiding. The notation is
from~\cite{Ma:bra}\cite{Ma:tra}\cite{Ma:bos} and
was recalled in Section~2. In our particular case in the category $\CM^H$ it
means
\alignn{rec}{ (h\tens b)(g\tens c)\equad &&= h\cdot g\bo \tens b\bo
c\CR(b\bt\tens g\bt)\nonumber \\
&&= h\t g\th \tens b\bo c \CR((Sh\o)h\th\tens  Sg\t)\CR(b\bt\tens
(Sg\o)g\fo) \\
\Delta (h\tens b)\equad &&= h\o\tens  b\Bo\bo\bo\tens h\t\bo \cdot
b\Bo\bt\tens b\Bt\CR(h\t\bt\tens b\Bo\bo\bt)\nonumber\\
&&= h\o\tens b\Bo\bo\tens h\t b\Bo\bt\tens b\Bt}
where we evaluated further in terms of $B,H$. The counit is the tensor product
one and there is an antipode as well. The coproduct comes out just the same as
for the bosonisation $H\rbiprod B$ (the usual cross coproduct by the coaction
of $H$ on $B$ as an object in $\CM^H$) because the transmutation procedure $B(\
,H)$ does not change the coalgebra. This is just the dual of the calculation of
$B\lcross B(H,H)$ in~\cite[Thm~3.2]{Ma:bos} for the module version.

These diagrammatic cross products and coproducts where algebraic
information `flows' along braids were introduced by the author in
\cite{Ma:bos}. Examples of braided module algebra structures are the
coregular representation which leads to braided-differentiation, and the
braided adjoint action~\cite[Prop.~3.1]{Ma:lie}
which leads to a theory of
braided Lie algebras. The comodule versions are the same with diagrams turned
up-side-down, e.g. the adjoint coaction~\cite[Appendix]{Ma:lin}, etc. We refer
to \cite{Ma:introp}\cite{Ma:introm} for
many basic results on braided (co)actions used in cross (co)products.

A trivial example of Proposition~6.1 is the case of
reconstruction of a super-Hopf algebra $\Z_2'\rcocross B$ from the category of
super $B$-modules and its forgetful functor. Here $\Z_2'$ is the dual of the
triangular Hopf algebra introduced by the author in
\cite[Prop~6.1]{Ma:exa}\cite[Ex.~1.1]{Ma:tra} as generator of the category
${\rm SuperVec}$ of superspaces with their $\Z_2$-graded transposition. This
 application  was generalised in
\cite{Ma:any} to generate the braided category of anyonic or $\Z_n$-graded
vector spaces introduced there. Unfortunately, in all these examples the
adjoint coaction of $H$ is trivial and $B(H,H)=H$ is viewed trivially in
$\CM^H$. Hence the algebra structure of $B(H,H)\rcocross B$ is the usual tensor
product one (and its cross coproduct the usual one as well). The result is a
braided group in $\CM^H$ just because $B$ is. The same applies for Hopf
algebras
in the braided category of $\Z$-graded vector spaces.

To give a more non-trivial example, let $q\in k^*$ and $H=GL_q(2)$ defined as
$k\<\alpha,\beta,\gamma,\delta,C^{-1}\>$ modulo the relations
\ceqn{qmat}{\alpha \beta=q^{-1}\beta \alpha,\quad \alpha \gamma=q^{-1}\gamma
\alpha,\quad \beta\delta=q^{-1}\delta \beta,\quad
\gamma\delta=q^{-1}\delta\gamma\\
\beta\gamma=\gamma\beta,\quad
\alpha\delta-\delta\alpha=(q^{-1}-q)\beta\gamma,\quad
C=\alpha\delta-q^{-1}\beta\gamma}
essentially as in~\cite{Dri}\cite{FRT:lie} for $SU_q(2)$. We equip it now with
dual quasitriangular structure determined by the associated solution $R$ of the
quantum Yang-Baxter equations. More precisely (for our application) we take $R$
with a non-standard normalisation as explained in\cite{Ma:lin}, fixed instead
by $\CR(C\tens C)=q^6$. Note therefore that one cannot set $C=1$ as one would
for the
usual dual quasitriangular Hopf algebra $SU_q(2)$.

The braided group $B(GL_q(2),GL_q(2))=BGL_q(2)$ is likewise a variant of the
braided group $BSU_q(2)$ introduced by the author in
\cite{Ma:eul}\cite{Ma:exa}. We define it as $k\<a,b,c,d,D^{-1}\>$ modulo the
relations
\ceqn{bmat}{ba=q^2ab,\quad ca=q^{-2}ac,\quad d a=ad,\qquad
bc=cb+(1-q^{-2})a(d-a)\\
d b=bd+(1-q^{-2})ab,\quad cd=d c+(1-q^{-2})ca,\quad D=ad-q^2cb}
It has `matrix' coproduct $\und\Delta\vecu=\vecu\tens\vecu$ and $\und\Delta
D=D\tens D$
when we regard the generators as a matrix $\vecu=\pmatrix{a&b\cr c&d}$. The
braided group antipode for $\vecu$ is as for $BSU_q(2)$ in
\cite{Ma:eul}\cite{Ma:exa} times $D^{-1}$. The braiding $\Psi$ between the
generators is also as listed for $BSU_q(2)$ in~\cite{Ma:eul}\cite{Ma:exa}. This
$BGL_q(2)$ lives in $\CM^{GL_q(2)}$ with a coaction which has the same `matrix
conjugation' form on the generators $\vecu$ as the right adjoint coaction of
$GL_q(2)$.

Let $B={\Bbb A}_q^2=k\<x,y\>/(yx-qxy)$ the $q$-deformed plane with right
coaction of $GL_q(2)$ given by transformation of the $(x,y)$ as a row vector by
the $GL_q(2)$ generators as a matrix, i.e. $\beta(x)=x\tens \alpha+y\tens
\gamma$ and $\beta(y)=x\tens\beta+y\tens\delta$. One of the first applications
of braided groups to physics was to show that this `quantum-braided plane'
${\Bbb A}_q^2$ is  a Hopf algebra in $\CM^{GL_q(2)}$ with linear `coaddition'
\cite{Ma:poi}
\ceqn{qplane}{ \Psi(x\tens x)=q^2 x\tens x,\quad \Psi(x\tens y)=q y\tens x\\
\Psi(y\tens y)=q^2 y\tens y,\quad \Psi(y\tens x)=q x\tens y+(q^2-1)y\tens x\\
\Delta x=x\tens 1+1\tens x,\quad \Delta y=y\tens 1+1\tens y,\quad \eps x=0=\eps
y,\\ Sx=-x,\quad Sy=-y.}
This result is due to the author in~\cite{Ma:poi}, where $GL_q(2)$ above is
formulated as $\widetilde{SU_q(2)}$, the `dilatonic' central extension.

We use the same matrix transformation for the braided coaction of $BGL_q(2)$ on
${\Bbb A}_q^2$. Under this, ${\Bbb A}_q^2$ becomes a right comodule algebra in
the braided category\cite[Prop.~3.7]{Ma:lin}.

\begin{example} The automorphism braided group $BGL_q(2)\rcocross {\Bbb A}_q^2$
in $\CM^{GL_q(2)}$ is generated by $BGL_q(2)$ and the quantum-braided plane
${\Bbb A}_q^2$ as subalgebras with the cross relations
\cmath{xa=ax,\quad ya=bx(q-q^{-1}) + ay,\quad xb=q^{-1}bx,\quad yb=qby,\quad
xc=qcx\\ yc=(1-q^{-2})(d-a)x + q^{-1}cy ,\quad xd=dx,\quad
yd=dy-q^{-2}(q-q^{-1})bx}
It has the matrix coproduct of $BGL_q(2)$ and
\[ \und\Delta x=x\tens a+y\tens c+1\tens x,\quad \und\Delta y=x\tens b+y\tens
d+1\tens
y\]
extended as a braided group in $\CM^{GL_q(2)}$.
\end{example}
\proof The cross relations are exactly the braided tensor product algebra as in
 (\ref{btens}), computed for the present setting  in terms of $R$ in
\cite[Lem.~3.4]{Ma:lin}. This gives the relations shown. For the coproduct we
know that we have the same form as the cross coproduct by the coaction of
$GL_q(2)$ on ${\Bbb A}_q^2$ but viewed now as a coaction of $BGL_q(2)$. To
extend the coproduct to products of the generators we use its
braided-multiplicativity, with $\Psi$ determined from the coaction. This was
computed in terms of $R$ in~\cite[Prop.~3.2]{Ma:lin} and in our case is
\ceqn{psiax}{\Psi(a\tens x)=x\tens a+(1-q^2)y\tens c,\  \Psi(b\tens
x)=q^{-1}x\tens b+(q-q^{-1})y\tens (a-d)\\ \Psi(c\tens x)=qx\tens c,\quad
\Psi(d\tens x)=x\tens d+ (1-q^{-2})y\tens c\\
\Psi(\pmatrix{a &b\cr c&d}\tens y)=y\tens \pmatrix{a&qb\cr q^{-1}c&d}}
while the braiding $\Psi(x\tens a)=a\tens x$ etc., has just the same form as
the cross relations already given. It is enough to specify the coproduct and
braiding on the generators since the braiding $\Psi$ itself extends
`multiplicatively' by functoriality and the Hexagon coherence identities, as
explained in~\cite{Ma:exa}. \endproof

The construction of linear braided groups such as the quantum-braided plane
works for general quantum planes associated to suitable matrix
data\cite{Ma:poi}. Another example is the 1-dimensional case $B={\Bbb
A}_q=k[x]$, the braided line\cite{Ma:sol}. Such `linear braided groups'  have
been very extensively studied since~\cite{Ma:poi} as the true foundation for
$q$-deformed geometry. See~\cite{Ma:varen} for a review. Their bosonisations
were used in~\cite{Ma:poi} to define inhomogeneous quantum groups and are also
extensively studied since then. Recently, the bosonisation construction has
been generalised so that both input and output are braided
groups\cite{Bes:cro}\cite{Dra:bos}.

\end{document}